\documentclass[11pt]{article}
\usepackage{amsfonts,amsmath,amssymb}
\newcommand{\eps}{\varepsilon}

\newcommand{\E}{\mathsf{E}}
\newtheorem{lemma}{Lemma}
\newcommand{\Fc}{\mathcal{F}}
\newtheorem{theorem}{Theorem}

\newtheorem{remark}{Remark}
\newcommand{\Pp}{\mathsf{P}}
\newcommand{\R}{\mathbb{R}}
\newcommand{\N}{\mathbb{N}}
\newcommand{\sgn}{\mathop{\rm sgn}}
\newcommand{\Nc}{\mathcal{N}}
\newcommand{\ONE}{{\bf 1}}
\newcommand{\Span}{\mathop{\rm span}}

\title{Exit asymptotics for small diffusion about an unstable equilibrium}
\author{Yuri Bakhtin\thanks{School of Mathematics, Georgia Institute of Technology, Atlanta, GA 30332, USA e-mail:{\it bakhtin@math.gatech.edu}}}

\begin{document}
\maketitle
\begin{abstract}
A dynamical system perturbed by white noise in a neighborhood of an unstable fixed point is considered. We obtain the exit asymptotics in the limit of vanishing noise intensity. This is a refinement of a result by Kifer (1981).
\end{abstract}

\section{Introduction}
Random perturbations of dynamical systems have been studied intensively for several decades,
see e.g. the classical book \cite{Freidlin--Wentzell-book}.
In particular, systems with
unstable equilibrium points including Hamiltonian and related flows have been considered, 
see e.g. recent works~\cite{Freidlin--Weber},\cite{Koralov}. See also~\cite{Armbruster--Stone--Kirk},~\cite{Ashwin--Borresen}, and~\cite{Afraimovich}  for results on noisy heteroclinic networks and their applications. 

The exit asymptotics for a neighborhood of an unstable fixed point
was studied in~\cite{Kifer}. It was shown that 
as the intensity $\eps$ of the white noise perturbation tends to $0$,
the exit distribution tends
to concentrate around the invariant manifold associated to the highest Lyapunov exponent $\lambda>0$, and
that the exit time $\tau$ is asymptotically equivalent to $\lambda^{-1}\ln(\eps^{-1})$ in probability. 

In this paper we prove a refinement of this asymptotics for additive isotropic noise. In particular, we show that 
$\tau - \lambda^{-1}\ln(\eps^{-1})$ converges almost surely to a random variable which we describe explicitly.

The approach we take also leads to
a simpler proof of the main theorem of~\cite{Kifer} for this setting. 
Our main result will be useful in analysis of vanishing noise asymptotics for dynamics with heteroclinic networks,
e.g. it provides a rigorous basis for some heuristic arguments from~\cite{Stone--Holmes}.

The paper is organized as follows. 
In Section~\ref{sec:main_result} we describe the setting and state our main result.
Its proof is given in Section~\ref{sec:main_proof} after a study of the linearized system
in Section~\ref{sec:linearization} and  estimates on closeness of the linear approximation to the original nonlinear system
in Section~\ref{sec:error_of_linearization}. Proofs of auxiliary lemmas are collected in Section~\ref{sec:aux_lemmas}.

{\bf Acknowledgments.} This work began after a talk on heteroclinic networks by Valentin Afraimovich. I am grateful to him for
several stimulating discussions. I would like to thank Vadim Kaloshin for a consultation on the Hadamard--Perron theorem.
I would also like to thank the referee whose comments and suggestions helped to improve the text a lot.
\section{The setting and the main result}\label{sec:main_result}

We suppose that there is a $C^2$-vector field $b:U\to\R^d$ defined on a bounded closed set
$U\subset \R^d$ equal to the closure of its own interior. 
This vector field generates a uniquely
defined flow $S^t$ associated with the ODE
\begin{equation*}
%\label{eq:ode}
\frac{d}{dt} S^tx=b(S^tx),\quad S^0x=x.
\end{equation*}
This flow is well-defined for all $t\in[0,T^U(x)]$ where $T^U(x)$ is the first time the solution hits 
$\partial U$:
\begin{equation*} 
T^U(x)=\inf\{t\ge0: S^tx\in\partial U\}.
\end{equation*}
A white noise perturbation of $S^t$ is given by the following SDE:
\begin{align} 
\label{eq:sde}
dX_\eps(t)&=b(X_\eps(t))dt+\eps dW(t),\\
X_\eps(0)&=x.\notag
\end{align}
Here $\eps>0$, and $W$ is a standard $d$-dimensional Wiener process defined on
a probability space $(\Omega,\Fc,\Pp)$. The SDE should be
understood in the integral sense:
\begin{equation}
\label{eq:sde_integral}
X_\eps(t)=x+\int_0^tb(X_\eps(t))dt+\eps W(t),
\end{equation}
and the (strong) solution can be obtained for $P$-almost every realization of~$W$.
We shall sometimes use the notation $S^{t}_{\eps,W} x$ to denote the solution
$X_\eps(t)$ of~\eqref{eq:sde_integral} to stress the dependence on the initial condition and the realization of the noise.
This solution is well-defined up to time $T^U_\eps(x)=T^U_{\eps,W}(x)$ which is a (random) stopping time defined as the first
time the solution hits~$\partial U$:
\begin{equation*} 
T^U_\eps(x)=T^U_{\eps,W}(x)=\inf\{t\ge 0: X_\eps(t)\in\partial U\}.
\end{equation*}
Let $G\subset U$ be a closed set with piecewise smooth boundary. For each $x\in G$ we can consider
equation~\eqref{eq:sde_integral} and define a stopping time
\begin{equation*}
\tau_\eps=\tau_\eps(x)=T^{G}_{\eps,W}(x)=\inf\{t\ge 0: X_\eps(t)\in\partial G\}
\end{equation*}
and the corresponding exit point
\begin{equation*}
H_\eps=H_\eps(x)=X_\eps(\tau_\eps).
\end{equation*}
With probability 1 we have $(\tau_\eps,H_\eps)\in\partial G\times \R_+$, and we
are going to study the asymptotics of this random vector in the limit as $\eps\to0$.

The limit behavior of $(\tau_\eps,H_\eps)$ depends very much on the vector field~$b$ and point $x$. 
We proceed to describe a setting which is slightly more restrictive than that of~\cite{Kifer}.

We shall assume that $0$ belongs to the interior of $G$, $b(0)=0$ and there are no other
equilibrium points in $G$.
We denote $A=J(0)$ where $J$ is the Jacobian matrix
\begin{equation*}
J=\left(\partial_i b_j\right)_{i,j=1,\ldots,d}.
\end{equation*} 
In this note we assume that $A$ has a simple positive eigenvalue $\lambda$  such that real parts
of all other eigenvalues are less than $\lambda$. We denote one of the two unit eigenvectors
associated to $\lambda$ by $v$. The vector subspace complement to~$v$ and spanned by all the other
Jordan basis vectors will be denoted by~$L$. Projections on $\Span\{v\}$ along $L$ and on $L$ along 
$\Span\{v\}$ will be denoted $\Pi_v$ and $\Pi_L$ respectively.

The Hadamard--Perron Theorem (see \cite[Theorem~6.2.8]{Katok--Hasselblatt} and \cite[Theorem~6.1]{Hartman})
implies that there is a locally $S^t$-invariant $C^1$-curve  $\gamma$ containing $0$ and tangent to $v$ at $0$.
%such that $\gamma(0)=0$ and $\gamma'(0)=v$. 
This curve is contained in the unstable manifold of the origin, and if all the other Lyapunov exponents are negative,
coincides with it. We shall assume that $\gamma\in C^2$ which is true
in many important cases.
In a small neighborhood of $0$, the curve~$\gamma$ can be represented as a graph of a map from $\Span\{v\}$ to $L$. 
For small $\delta$ we shall denote by $\gamma(\delta)$ the point $x$ on $\gamma$ such that $\Pi_v x=\delta v$.   
Notice that $|\Pi_L\gamma(\delta)|=O(\delta^2)$ as $\delta\to0$, where $|\cdot|$ denotes the Euclidean norm in $\R^d$.

We also assume that $\gamma$ intersects $\partial G$ transversally at two points $q_-$ and $q_+$ so that
the part of $\gamma$ connecting $q_-$ and $q_+$ does not
intersect $\partial G$ and contains points $\gamma(-\delta),0,\gamma(\delta)$ (in this order) for some $\delta>0$. 

We shall need the quantities $h_+$ and $h_-$ defined via:
\begin{equation}
\label{eq:h}
h_{\pm}=\lim_{\delta\to 0} \left(\frac{\ln\delta}{\lambda}+t(\pm\delta,q_{\pm})\right),
\end{equation}
where $t(\delta,q_+)$ and $t(-\delta,q_-)$ denote the times to get from $\gamma(\delta)$ to
$q_+$ and from $\gamma(-\delta)$ to
$q_-$ respectively: 
\begin{equation}
t(\pm\delta,q_\pm)=T^G(\gamma(\pm\delta)),
\label{eq:definition_of_t_q_+}
\end{equation}
so that
$
S^{t(\pm\delta,q_\pm)} \gamma(\pm\delta)=q_{\pm}.
$

\begin{lemma}\label{lm:h_well_defined} The numbers $h_\pm$ are well-defined by~\eqref{eq:h}, i.e. finite limits in the r.h.s. exist.
\end{lemma}

A proof of this lemma is given in Section~\ref{sec:aux_lemmas}, and we proceed now to our main result.
\begin{theorem}\label{th:main} Suppose $x$ belongs to the exponentially stable manifold of $0$, i.e. there are positive constants $C$ and $\mu$
such that
\begin{equation}
\label{eq:stable_manifold}
|S^tx|\le Ce^{-\mu t},\quad t\ge0.
\end{equation} 
Then the following holds:
\begin{enumerate}
\item There is a positive number $\sigma=\sigma(x)$ and a standard Gaussian random variable~$\Nc$ defined on the
probability space $(\Omega,\Fc,\Pp)$ such that with probability~$1$
\begin{equation*}
H_\eps\to q_+\ONE_{\{\Nc>0\}}+q_-\ONE_{\{\Nc<0\}},
\end{equation*}
and
\begin{equation*}
\tau_\eps-\frac{\ln(1/\eps)}{\lambda}\to h_+\ONE_{\{\Nc>0\}}+h_-\ONE_{\{\Nc<0\}}-\frac{\ln(\sigma|\Nc|)}{\lambda}.
\end{equation*}
\item The distribution of the random vector $(H_\eps, \tau_\eps-\frac{\ln(1/\eps)}{\lambda})$
converges weakly to
\begin{equation*}
\frac{1}{2}\delta_{q_+}\times \mu_{h_+,\sigma} +\frac{1}{2}\delta_{q_-}\times\mu_{h_-,\sigma},
\end{equation*}
where $\mu_{h,\sigma}$ is the distribution of 
\begin{equation*}
h-\frac{\ln(\sigma|\Nc|)}{\lambda}.
\end{equation*}
\item If $x=0$, then $\sigma=(2\lambda)^{-1/2}$.
\end{enumerate}
\end{theorem}
%\begin{remark} An explicit expression for $\sigma$ can also be obtained for the generic case \end{remark}
\begin{remark} Notice that if $A$ has no negative eigenvalues, then the only point~$x$ satisfying \eqref{eq:stable_manifold}
for some $C,\mu>0$ is the origin, i.e. the stable manifold is trivial, and our theorem is applicable only for the diffusion started
at $x=0$. In the opposite case, where there is at least
one negative eigenvalue, the Hadamard--Perron theorem mentioned above guarantees the existence of a nontrivial stable manifold which plays the role of the
unstable one for the system in the reverse time. Notice also that in the latter situation one can choose $\mu$ to be a constant independent of $x$ (namely,
take any negative number that is closer to zero than any negative Lyapunov exponent), but $C$ depends on $x$ essentially.
\end{remark}

\begin{remark} It is an interesting phenomenon that in the situation where there is a variety of unstable
directions, the system chooses the most 
unstable one so that the limiting dynamics is practically
1-dimensional if the leading eigenvalue of the linearization is real and simple. This was observed 
in~\cite{Kifer}, where, in fact, a more general situation was considered as well. The leading eigenvalue
$\lambda$ was not necessarily assumed real and simple. We can easily extend our approach to
recover the main result of~\cite{Kifer}: 
the exit time grows as $\lambda^{-1}\ln(\eps^{-1})$ and the exit
measure tends to concentrate at the intersection of $\partial G$ and the invariant manifold
corresponding to $\lambda$. However, without our assumptions on $\lambda$, the exit asymptotics analogous to
Theorem~\ref{th:main} is more
complicated and depends, in particular, on the shape of the set $G$.

\begin{remark} The random variable $\Nc$ is constructed explicitly in the proof of the theorem.
\end{remark}

\end{remark}

\section{Linearization}\label{sec:linearization}
We start our study of the SDE~\eqref{eq:sde} with the analysis of its linearization:
\begin{equation}
\label{eq:X_tilde}
\tilde X_\eps(t)=S^tx+\eps Y(t),
\end{equation}
where $Y$ solves the equation in variations:
\begin{align*}
dY(t)&=A(t)Y(t)dt+dW(t),\\
Y(0)&=0.
\end{align*}
Here $A(t)=J(S^tx)$. 

The main result of this section is the following lemma.
\begin{lemma}\label{lm:convergence_to_gaussian}
There is a centered nondegenerate Gaussian random variable $N$, an a.s.-finite random variable $C_1$
and a number $\rho>0$ such that with probability 1, for all $t$,
\begin{equation*}
|e^{-\lambda t} Y(t)- Nv|\le C_1e^{-\rho t}.
\end{equation*}
\end{lemma}
\begin{remark} The Gaussian random variable $N$ will be used to construct  $\sigma$ and~$\Nc$ that appear in the statement of Theorem~\ref{th:main}.
Namely, $N=\sigma\Nc$, see Section~\ref{sec:main_proof}.

\end{remark}
Proof. Let $Z$ be the solution of
\begin{align*}
dZ(t)&=AZ(t)dt+dW(t),\\
Z(0)&=0.
\end{align*}
Then
\begin{equation*}
Z(t)=\int_0^t e^{A(t-s)}dW(s),
\end{equation*}
see~\cite[Section 5.6]{Karatzas--Shreve} for a treatment of stochastic linear equations. 
Let us denote $V(t)=e^{-\lambda t} Z(t)$.
Since $e^{A(t-s)}v=e^{\lambda (t-s)} v$, and 
\begin{equation}
\label{eq:exponential_action_on_L}
|e^{A(t-s)}u|<C_2e^{(\lambda-\Delta)(t-s)}|u|
\end{equation} 
 for some positive constants $\Delta,C_2$ and all $u\in L$,
we have
\begin{equation*}
\label{eq:asymptotics_with_constant_matrix}
V(t)\to\Pi_{v} \int_0^\infty e^{-\lambda s} d W(s),\quad\mbox{\rm as}\ t\to\infty,
\end{equation*}
%where
%\begin{equation*}
%V(t)=e^{-\lambda t} Z(t),
%\end{equation*}
and the convergence is exponentially fast.

Now let $D(t)=Y(t)-Z(t)$. Then
\begin{equation*}
\frac{d}{dt}D(t)=AD(t)+(A(t)-A)Y(t),
\end{equation*}
so that
\begin{align*}
D(t)=&\int_0^te^{A(t-s)}(A(s)-A)Y(s)ds\\
    =&\int_0^te^{\lambda(t-s)}\Pi_v (A(s)-A)Y(s)ds
     + \int_0^te^{A(t-s)}\Pi_L(A(s)-A)Y(s)ds.
\end{align*}
This implies
\begin{multline}
\label{eq:prelimit_for_D(t)}
e^{-\lambda t}D(t)= \int_0^t e^{-\lambda s}\Pi_v(A(s)-A)Y(s)ds\\
     +e^{-\lambda t}  \int_0^te^{A(t-s)}\Pi_L(A(s)-A)Y(s)ds.
\end{multline}

To estimate the r.h.s. we write
\begin{equation}
\label{eq:convolution_for_Y}
Y(s)=\int_0^s\Phi_r(s)dW(r),
\end{equation}
where $\Phi_r(s)$ denotes the 
fundamental matrix solving the homogeneous system
\begin{align}
\frac{d}{ds}\Phi_{r}(s)&=A(s)\Phi_{r}(s)\label{eq:Phi1},\\
\Phi_r(r)&=I.\label{eq:Phi2}
\end{align}

For a matrix $B$, we denote $|B|=\sup_{|x|\le1}|Bx|$.
\begin{lemma}\label{lm:growth_of_Phi} 
For any $\alpha>0$ there is a constant $K_\alpha$ such that 
\begin{equation*}
|\Phi_s(t)|\le K_\alpha e^{(\lambda+\alpha)(t-s)}
\end{equation*}
for all $t,s$ with $t>s>0$.
\end{lemma}

We prove this lemma in Section~\ref{sec:aux_lemmas}. An almost immediate implication is the following statement:
%This also implies that 
%\begin{equation*}
%\limsup_{s\to\infty}\frac {\ln|Y(s)|}{s}\le\lambda,
%\end{equation*}
%see [???], so 
\begin{lemma}\label{lm:growth_of_Y}
For any $\alpha>0$ there is an a.s.-finite  random constant $\tilde K_\alpha$ such that with probability 1,
\begin{equation*}
%\label{eq:growth_of_Y}
|Y(s)|\le \tilde K_\alpha e^{(\lambda+\alpha)s}.
\end{equation*}
for all $s\ge 0$. 
\end{lemma}
The proof of this lemma is also given in Section~\ref{sec:aux_lemmas}.
It is important that the positive number $\alpha$ can be chosen arbitrarily small. In fact, we shall use Lemmas~\ref{lm:growth_of_Phi} and \ref{lm:growth_of_Y} for
 $\alpha<\mu$.

Since $x$ belongs to the stable manifold of the origin, we have
\begin{equation*}
|A(s)-A|\le C_3 e^{-\mu s}
%\label{eq:rate_of_approximation_for_A}
\end{equation*}
for some $C_3$ and all $s\ge 0$, where $\mu$ was introduced in~\eqref{eq:stable_manifold}. Therefore, Lemma~\ref{lm:growth_of_Y} implies that
as $t\to\infty$, the first integral in~\eqref{eq:prelimit_for_D(t)}
exponentially converges to
\begin{equation*}
\Pi_v \int_0^\infty e^{-\lambda s}(A(s)-A)Y(s)ds.
\end{equation*} 
The same considerations and~\eqref{eq:exponential_action_on_L} imply that
the second integral in~\eqref{eq:prelimit_for_D(t)} converges to $0$
exponentially fast.

Therefore,
\begin{equation}
\label{eq:explicit_expression_for_N}
\lim_{t\to\infty} e^{-\lambda t} Y(t)\stackrel{a.s.}{=}\Pi_{v} \left[\int_0^\infty e^{-\lambda s} d W(s)+\int_0^\infty e^{-\lambda s}(A(s)-A)Y(s)ds\right].
\end{equation}

The r.h.s. is a Gaussian random variable with distribution concentrated on $\Span\{v\}$ since it is a finite linear functional 
of the Wiener process $W$. Our proof will be complete as soon as we show that this linear functional is non-degenerate.
Using~\eqref{eq:convolution_for_Y} we rewrite the r.h.s. of \eqref{eq:explicit_expression_for_N} as
\begin{align}
\Pi_{v} \left[\int_0^\infty e^{-\lambda r} d W(r)+\int_0^\infty e^{-\lambda s}(A(s)-A)\int_0^s\Phi_r(s)dW(r)ds\right]\notag\\
=\Pi_{v} \left[\int_0^\infty \left(e^{-\lambda r}I +\int_r^\infty e^{-\lambda s}(A(s)-A)\Phi_r(s)ds\right)\, dW(r)\right],
\label{eq:nondegenerate_gaussian}
\end{align}
where $I$ denotes the unit matrix.

Let us take a positive $\alpha<\mu$. 
Lemma~\ref{lm:growth_of_Phi}  implies that
\begin{align*}
\left|\int_r^\infty e^{-\lambda s}(A(s)-A)\Phi_r(s)ds\right|&\le C_3K_\alpha\int_r^\infty e^{-\lambda s}e^{-\mu s}e^{(\lambda+\alpha)(s-r)}ds\\&\le\frac{C_3K_\alpha}{\mu-\alpha}e^{-(\lambda+\mu)r}, 
\end{align*}
and the expression in the stochastic integral in the r.h.s. of~\eqref{eq:nondegenerate_gaussian} cannot be identically equal to zero which completes the proof of Lemma~\ref{lm:convergence_to_gaussian}.

\begin{lemma}\label{lm:variance_for_zero_initial_data} If $x=0$, then $\E N^2=1/(2\lambda)$, where $N$ is the centered Gaussian random variable defined
in Lemma~\ref{lm:convergence_to_gaussian}.
\end{lemma}

Proof. If $x=0$, then $A(t)=A$ for all $t\ge0$. Therefore, the second term in the r.h.s. of~\eqref{eq:explicit_expression_for_N}
vanishes, and the variance of the first term equals $\int_0^\infty e^{-2\lambda s}ds=1/(2\lambda)$ due to It\^o's isometry, 
see \cite[Lemma 3.5]{Oksendal}.

For every $\delta>0$ we shall need a stopping time 
\begin{align*}
\tau(\tilde X_\eps,\delta,v)&=\inf\{t>0:\ |\Pi_{v}(\tilde X_\eps(t)-S^tx)|\ge\delta \}\\&=\inf\{t>0:\ \eps|\Pi_{v}Y(t)|\ge\delta \},
\end{align*}
where $\tilde X_\eps$ is defined in~\eqref{eq:X_tilde}.

\begin{lemma}\label{lm:exit_for_linear} For any $\delta>0$,
\begin{equation*}
\lim_{\eps\to 0}\left[\tau(\tilde X_\eps,\delta,v)-\frac{\ln\left(\frac{\delta}{\eps|N|}\right)}{\lambda}\right]\stackrel{a.s.}{=}0,
\end{equation*}
where $N$ is the centered Gaussian random variable defined
in Lemma~\ref{lm:convergence_to_gaussian}.
\end{lemma}

Proof. Obviously, $\tau=\tau(\tilde X_\eps,\delta,v)\stackrel{a.s.}{\to}\infty$ as $\eps\to 0$, so that
Lemma~\ref{lm:convergence_to_gaussian} implies
\begin{equation*}
\delta = \eps|\Pi_v  Y(\tau)|\sim \eps e^{\lambda \tau} |N|,
\end{equation*}
and the claim follows.

\begin{lemma}\label{lm:linear_exit_2} There is a positive number $\beta$ such that for any~$\delta>0$ there is an a.s.-finite random variable $C_4=C_4(\delta)$ such that
%\begin{equation}
%\lim_{\eps\to0}\Pp\{\eps\Pi_{v}Y(\tau)=\pm\delta v\}=\frac{1}{2},
%\label{eq:symmetry}
%\end{equation}
%and 
with probability 1 
\begin{equation*}
%\label{eq:eps-closeness}
\limsup_{\eps\to 0} \frac{|\eps\Pi_{L}Y(\tau)|}{\eps^{\beta}}\le C_4.
\end{equation*}
\end{lemma}

Proof. 
%Equation~\eqref{eq:symmetry} follows from Lemma~\ref{lm:convergence_to_gaussian} and
%the symmetry of the Gaussian distribution. Next, 
Lemmas~\ref{lm:convergence_to_gaussian} and~\ref{lm:exit_for_linear} imply that
\begin{equation*}
|\eps\Pi_{L}Y(\tau)|\le C_1\eps e^{\lambda\tau}e^{-\rho\tau}\sim C_1\frac{\delta}{|N|}\cdot\left(\frac{\eps|N|}{\delta}\right)^{\rho/\lambda},
\end{equation*}
which proves our claim with $\beta=\rho/\lambda$ and $C_4(\delta)=C_1(\delta/|N|)^{1-\rho/\lambda}$.

\section{The error of the linear approximation}\label{sec:error_of_linearization}
In this section, we are going to compare the nonlinear diffusion process
$X_\eps$ to its Gaussian linearization $\tilde X_\eps$ considered in Section~\ref{sec:linearization}.
%Let the stopping time 
%$\tau(\tilde X_\eps,\delta,v)$ be defined as in the last section.

\begin{lemma}
\label{lm:delta_squared_error_multidim} There is a number $\delta_0>0$ and a random variable $C_5>0$
such that, with probability 1, if $\delta\in(0,\delta_0)$, then
\begin{equation*}
\limsup_{\eps\to 0}|X_\eps(\tau(\tilde X_\eps,\delta,v)) - \tilde X_\eps(\tau(\tilde X_\eps,\delta,v))|< C_5\delta^2.
\end{equation*}
\end{lemma}

Proof.
%We shall often use $\tau$ instead of $\tau(Z_\eps,\delta,v)$ for brevity.
In differential notation, the evolution of $\tilde X_\eps$ is given by
\begin{align*}
d\tilde X_\eps(t)&=b(S^tx)dt+\eps dY(t)\\&=b(S^tx)dt+\eps A(t)Y(t)dt+\eps dW(t).
\end{align*}
Using $Y(t)=(\tilde X_\eps(t)-S^tx)/\eps$, we obtain
\begin{equation*}
d\tilde X_\eps(t)=b(S^tx)dt+A(t)(\tilde X_\eps(t)-S^tx)dt+\eps dW(t).
\end{equation*}
Let us introduce $U_\eps(t)=X_\eps(t)-\tilde X_\eps(t)$, so that $U_\eps(0)=0$ and
\begin{equation*}
\frac{d}{dt} U_\eps(t)=b(X_\eps(t))-b(S^tx)-A(t)(\tilde X_\eps(t)-S^tx).
\end{equation*}
Since $b\in C^2$, we have
\begin{equation*}
b(X_\eps(t))=b(S^tx)+A(t)(X_\eps(t)-S^tx)+Q(S^tx, X_\eps(t)-S^tx)
\end{equation*}
where
\begin{equation*}
|Q(y, z)|\le C_6 |z|^2
\end{equation*}
for a constant $C_6$ and all $y,z$, so that
\begin{equation*}
\frac{d}{dt} U_\eps(t)=A(t)U_\eps(t)+Q(S^tx, X_\eps(t)-S^tx).
\end{equation*}

Variation of constants yields:
\begin{equation*}
V(t)=\int_0^t \Phi_s(t)Q(S^sx, X_\eps(s)-S^sx)ds,
\end{equation*}
where $\Phi_s(t)$ is defined in~\eqref{eq:Phi1}--\eqref{eq:Phi2}.
Since 
\begin{equation*}
|Q(S^sx, X_\eps(s)-S^sx)|\le C_6 |U_\eps(s)+ \tilde X_\eps(s)-S^sx|^2\le 2C_6|U_\eps(s)|^2+2C_6\eps^2|Y(s)|^2,
\end{equation*}
Lemma~\ref{lm:growth_of_Phi} implies that for any $\alpha>0$,
\begin{equation*}
|U_\eps(t)|\le 2K_\alpha C_6 \int_0^t e^{(\lambda+\alpha)(t-s)} (|U_\eps(s)|^2+\eps^2|Y(s)|^2)ds,
\end{equation*}
so that $|U_\eps(t)|\le u_\eps(t)$, where $u_\eps$ solves
\begin{align}
\label{eq:bounding_solution_u}
\frac{d}{dt} u_\eps(t)&=(\lambda+\alpha)u_\eps(t)+2K_\alpha C_6 u_\eps^2(t)+2K_\alpha C_6\eps^2|Y(t)|^2, \\
u_\eps(0)&=0.\notag
\end{align}
Obviously, $u_\eps$ is a monotone nondecreasing
function. Let us choose $\alpha<\lambda/2$ and denote 
\begin{equation*}
c=\frac{\frac{1}{2}\lambda-\alpha}{2K_\alpha C_6}.
\end{equation*}
If $|u_\eps(t)|\le c$, then
\begin{equation*}
(\lambda+\alpha)u_\eps(t)+2K_\alpha C_6 u_\eps^2(t)\le \frac{3}{2}\lambda u_\eps(t). 
\end{equation*}
Therefore,
\begin{equation*}
\ONE_{\{u_\eps(t)\le c\}}\frac{d}{dt} u_\eps(t)\le \frac32\lambda u_\eps(t) +2K_\alpha C_6\eps^2|Y(t)|^2,
\end{equation*}
so that
\begin{equation*}
u_\eps(t)\ONE_{\{u_\eps(t)\le c\}}\le 2K_\alpha C_6\eps^2 \int_0^t
e^{\frac32\lambda(t-s)} |Y_\eps(s)|^2 ds.
\end{equation*}
Since $|Y_\eps(t)|\sim e^{\lambda t}|N|$, the r.h.s.
is asymptotically equivalent to
\begin{equation*}
2K_\alpha C_6 N^2\eps^2\int_0^t e^{\frac32\lambda(t-s)}e^{2\lambda s}ds
\sim 2K_\alpha C_6 N^2\eps^2 e^{2\lambda t},
\end{equation*}
which implies that
\begin{equation}
\limsup_{t\to\infty}\frac{u_\eps(t)\ONE_{\{u_\eps(t)\le c\}}}{2K_\alpha C_6N^2\eps^2 e^{2\lambda t}}\le 1
\label{eq:localization_u_eps}
\end{equation}
uniformly in $\eps>0$.
Next, let us consider $\tau(u_\eps,c)=\inf\{t\ge0: u_\eps(t)\ge c\}$. Monotonicity of the r.h.s of~\eqref{eq:bounding_solution_u} in $\eps$ implies
$\tau(u_\eps,c)\to\infty$ as $\eps\to 0$.
Since $|u_\eps(s)|\le c$ for all $s\le \tau(u_\eps,c)$ and
$|u_\eps(\tau(u_\eps,c))|=c$,
formula~\eqref{eq:localization_u_eps} implies
\begin{equation*}
\limsup_{\eps\to 0}\frac{c }{2K_\alpha C_6 N^2\eps^2 e^{2\lambda \tau(u_\eps,c)}}\le 1,
\end{equation*}
i.e.
\begin{equation*}
\liminf_{\eps\to0} \left[\tau(u_\eps,c)-\frac{\ln(\frac{1}{|N|\eps})}{\lambda}-\frac{\ln(\frac{c}{2K_\alpha C_6})}{2\lambda}\right]\ge 0.
\end{equation*}
%where $f(\eps)\gtrapprox g(\eps)$ means $\liminf_{\eps\to0}(f(\eps)-g(\eps))\ge 0$.
The last relation and Lemma~\ref{lm:exit_for_linear} imply that for sufficiently small $\delta_0$ and all $\delta\in(0,\delta_0)$, there is an
$\eps_0=\eps_0(\delta)$ such that
if  $0<\eps<\eps_0$, then
\begin{equation*}
\tau(\tilde X_\eps,\delta,v)< \tau(u_\eps,c).
\end{equation*}
Now \eqref{eq:localization_u_eps} implies that for these
$\delta$ and sufficiently small $\eps$
\begin{equation*}
u_\eps(\tau(\tilde X_\eps,\delta,v))\le 3K_\alpha C_6 N^2\eps^2 e^{2\lambda \tau(\tilde X_\eps,\delta,v)}\le 4K_\alpha C_6\delta^2,
\end{equation*}
where in the last inequality we used~Lemma~\ref{lm:exit_for_linear} again. The proof is complete.

%We shall now introduce a new pair of an exit time and exit point:
%\begin{align*}
%\tau(X_\eps,\delta,v)&=\inf\{t>0:\ |\Pi_{v}(X_\eps(t)-S_tx)|\ge\delta \},\\
%H(X_\eps,\delta,v)&=X_\eps(\tau(X_\eps,\delta,v)).
%\end{align*}
%\begin{lemma}
%\label{th:delta_squared_error_multidim} With probability 1 there are constants $C,\delta_0>0$ and a function
%$\eps_0:(0,\delta_0)\to(0,\infty)$
%such that if $\delta\in(0,\delta_0)$ and $\eps\in(0,\eps_0(\delta))$, then
%\begin{equation*}
%|Y_\eps(\tau(Z_\eps,\delta,v)) - Z_\eps(\tau(Z_\eps,\delta,v))|\le C\delta^2.
%\end{equation*}
%\end{lemma}

\section{Proof of the main result}\label{sec:main_proof}
We begin with auxiliary well-known statements. 
The first statement estimates closeness of perturbed trajectories to the orbits of the unperturbed system. 
\begin{lemma}\label{lm:closeness_to_unperturbed} Let $W^*(t)=\sup_{s\in[0,t]}|W(s)|$. Then, for any $y\in U$, any $t<T^U(y)$, for a.e. Wiener trajectory $W$ and
$\eps<\eps_0(W)$,
\begin{equation*}
|S^t_{\eps,W}y-S^ty|\le \eps W^*(t) e^{Mt},
\end{equation*}
where $M$ is the Lipschitz constant of $b$ on $U$.
\end{lemma}
Proof of Lemma~\ref{lm:closeness_to_unperturbed}. Denote $V_\eps(t)=S^t_{\eps,W} y-S^tx$. 
Then
\begin{equation*}
V(t)=\int_0^t (b(S^s_{\eps,W} y)-b(S^sy))ds+\eps W(s),
\end{equation*}
so that
\begin{equation*}
|V_\eps(t)|\le\int_0^t M|V_\eps(s)|ds+\eps |W(s)|,
\end{equation*}
and the lemma follows from Gronwall's inequality and a simple localization argument. 

The next statement will estimate the closeness to the invariant curve~$\gamma$. 
We need more notation.
%For all points $y\in G$ that are sufficiently close to $\gamma$ and do not belong to the stable manifold of the origin,
%$T^{G}(y)$ is a well-defined finite positive number, and we can define the exit point as
%\begin{equation*} 
%H_0(y)=S^{T^G(y)}y.
%\end{equation*}
%In particular, for small $\delta>0$
%\begin{equation*} 
%H_0(\gamma(\pm\delta))=q_{\pm}.
%\end{equation*}

For $K>0$ we introduce two sets
\begin{equation*}
D^{\pm}(\delta,K)=\{x\in\R^d:|x\mp\delta v|\le K\delta^2\}.
\end{equation*}

We shall need a closed set $G'\subset U$ with the following properties: the boundary of $G'$ is piecewise smooth, $G$ is contained in the interior of $G'$, and
$\gamma$ intersects $\partial G'$ transversally.

\begin{lemma}\label{lm:exit_times_G_prime}
For any $K>0$ and sufficiently small $\delta$, there are positive numbers $T^+=T^+(\delta)$ 
and $T^-=T^-(\delta)$ such that
\begin{equation*}
S^{T^\pm} D^{\pm}(\delta,K)\subset U\setminus G'.
\end{equation*}
For any $K>0$ 
\begin{equation*}
\lim_{\delta\to 0}\sup_{t\le T^{G'}(\gamma(\pm\delta))}\,\sup_{y\in D^{\pm}(\delta,K)}|S^t y-S^t(\gamma(\pm\delta))|=0.
\end{equation*}
\end{lemma}
This lemma follows from the graph transform method of
constructing the invariant manifold $\gamma$ (see e.g. a version of Hadamard--Perron Theorem
and its proof in
\cite[Section 6.2]{Katok--Hasselblatt}).

%The following is a straightforward consequence of the last lemma:
%\begin{lemma}
%\label{lm:closeness_to_unstable_manifold} 
%For any $K>0$, 
%\begin{equation*}
%\lim_{\delta\to 0}\,\sup_{y\in D^{\pm}(\delta,K)} |H_0(y)-H_0(\gamma(\pm\delta))|= 0,
%\end{equation*}
%and
%\begin{equation*}
%\lim_{\delta\to 0}\,\sup_{y\in D^{\pm}(\delta,K)} |T^G(y)-T^G(\gamma(\pm\delta))|= 0.
%\end{equation*}
%\end{lemma}

Proof of Theorem~\ref{th:main}. For any $y\in G$ and any time $\nu\ge 0$ we define $H_\eps^{\nu}(y)$
and $\tau_\eps^\nu(y)$ analogously to  $H_\eps(y)$
and $\tau_\eps(y)$, but using shifted trajectories $W(\nu+\cdot)-W(\nu)$ instead of $W(\cdot)$.

For sufficiently small $\delta,\eps>0$,
\begin{align}
H_\eps(x)&=H_\eps^{\tau(\tilde X_\eps,\delta,v)}(X_\eps(\tau(\tilde X_\eps,\delta,v))), \label{eq:H_decomposition}\\
\tau_\eps(x)&=\tau(\tilde X_\eps,\delta,v)+\tau_\eps^{\tau(\tilde X_\eps,\delta,v)}(X_\eps(\tau(\tilde X_\eps,\delta,v))).
\label{eq:tau_decomposition}
\end{align}

Let us now define the nondegenerate Gaussian random variable $N$
according to Lemma~\ref{lm:convergence_to_gaussian}.
Lemma~\ref{lm:linear_exit_2} (estimating the closeness of the linearized
process to $\sgn(N)\delta v$ at the exit time $\tau(\tilde X_\eps,\delta,v)$) and Lemma~\ref{lm:delta_squared_error_multidim} (estimating the 
closeness of the nonlinear process to the linearized process
at the time $\tau(\tilde X_\eps,\delta,v)$) imply that there is a constant $\delta_0>0$ and a positive a.s.-finite random variable~$C_7$
such that with probability 1 for $\delta\in(0,\delta_0)$  we have 
\begin{equation}
\label{eq:diff-with-sgn(N)delta.v}
\limsup_{\eps\to 0}\left| X_\eps(\tau(\tilde X_\eps,\delta,v))-\sgn(N) \delta v\right|<C_7\delta^2.
\end{equation}

To estimate the effect of the noise for the evolution along $\gamma$, we take a~$K>0$ and write
\begin{align}
\sup_{y\in D^{\pm}(\delta,K)}|S^t_{\eps,W}y &- S^t(\gamma(\pm\delta))|\notag \\&\le 
\sup_{y\in D^{\pm}(\delta,K)}|S^t_{\eps,W}y-S^t y|+
\sup_{y\in D^{\pm}(\delta,K)}|S^t y-S^t(\gamma(\pm\delta))|\notag \\
&\le \eps W^*(t)e^{Mt}+r(\pm\delta,K)\label{eq:eps_and_delta_estimate}
\end{align}
if all the involved processes stay within $G'$ up to time $t$.
Here we used Lemma~\ref{lm:closeness_to_unperturbed} to bound the first
term and denoted by $r(\pm\delta,K)$ the second term. Notice that
for each $K>0$ we have $r(\pm\delta,K)\to0$ as $\delta\to0$ due to Lemma~\ref{lm:exit_times_G_prime}.

So, we need an estimate on the exit time from $G'$.
Let us fix $K>0$ and small $\delta>0$.
Due to the continuous dependence of $S^t_{\eps,W}y$ on $y$ and~$\eps$, Lemma~\ref{lm:exit_times_G_prime}
allows us to choose constant times $\tilde T(\pm\delta,K)>0$ and $\eps_0=\eps_0(W)>0$ such that
for all $\eps\in(0,\eps_0)$, we have
$(S^{\tilde T(\pm\delta,K)}_{\eps,W} D(\pm\delta,K))\cap G'=\emptyset$. 
In particular, for fixed $\delta$ and $K$, and for almost every Wiener trajectory $W$, we have
\begin{equation*}
\label{eq:finite_horizon_for_tau}
\limsup_{\eps\to0} \sup_{y\in{D^\pm(\delta,K)}} T^G_{\eps,W}(y)\le \tilde T(\delta,K),
\end{equation*}
and \eqref{eq:eps_and_delta_estimate} implies that
\begin{equation}
\limsup_{\eps\to0}\sup_{y\in D^{\pm}(\delta,K)}|S^t_{\eps,W}y - S^t(\gamma(\pm\delta))|\le r(\pm\delta,K)
\label{eq:W-closeness}
\end{equation}
for all $t\le \tilde T(\delta,K)$. 

Since for any $y\in D^\pm(\delta,K)$, almost every $W$ and all sufficiently small~$\eps>0$, we have  $T^G_{\eps,W}(y)\le \tilde T(\delta,K)$,
we can combine \eqref{eq:W-closeness} with the strong Markov property and estimate
\eqref{eq:diff-with-sgn(N)delta.v} to see that 
\begin{equation*}
\limsup_{\eps\to 0}\bigl|H_\eps^{\tau(\tilde X_\eps,\delta,v)}(X_\eps(\tau(\tilde X_\eps,\delta,v)))-H_0(\gamma(\sgn N\delta))\bigr|\ONE_{\{C_7< K\}} \le r_1(\pm\delta,K),
%\label{eq:conv_of_H}
\end{equation*}
and
\begin{equation*}
\limsup_{\eps\to 0}|\tau_\eps^{\tau(\tilde X_\eps,\delta,v)}(X_\eps(\tau(\tilde X_\eps,\delta,v)))-T^G(\gamma(\sgn N\delta))|\ONE_{\{C_7< K\}}\le r_1(\delta,K),
%\label{eq:conv_of_tau}
\end{equation*}
for a deterministic function $r_1(\delta,K)>0$ such that $\lim_{\delta\to 0}r_1(\delta,K)=0$ for any fixed $K>0$.

Therefore, Lemma~\ref{lm:exit_for_linear}  and
identities \eqref{eq:H_decomposition} and \eqref{eq:tau_decomposition} imply that
for fixed $\delta$ and $K$, on $\{C_7<K\}$ we have 
\begin{equation*}
\limsup_{\eps\to 0} \left|\tau_\eps(x)-\frac{\ln\left(\frac{\delta}{\eps|N|}\right)}{\lambda}-T^G(\gamma(\sgn N\delta))\right|\le r_2(\delta,K),
\end{equation*}
and
\begin{equation*}
\limsup_{\eps\to 0} |H_\eps(x)-H_0(\gamma(\sgn N\delta))|\le r_2(\delta,K)
\end{equation*}
with $\lim_{\delta\to 0} r_2(\delta,K)=0$. 

Since $\{C_7<K\}\uparrow \Omega$ as $K\to\infty$, Part~1 of Theorem~\ref{th:main} follows with $\sigma=\sqrt{\E N^2}>0$ and 
$\Nc=N/\sigma$. 

Part~2 follows from Part~1 as soon as we notice that $|N|$ and $\sgn N$ are independent, the latter taking values $\pm1$ with probability $1/2$.

Part~3 of Theorem~\ref{th:main} follows from Lemma~\ref{lm:variance_for_zero_initial_data},
and the proof is complete.

\section{Auxiliary Lemmas}\label{sec:aux_lemmas}
Proof of Lemma~\ref{lm:h_well_defined}.
Let us prove that $h_+$ is well-defined by \eqref{eq:h} (the same proof works for $h_-$ as well). Let $z(t)=|\Pi_v S^{-t}q_+|$. 
There is $t_0>0$ such that on the semiline $(t_0,+\infty)$ the function
$z(t)$ is monotone decreasing and satisfies
\begin{equation*}
\dot z(t) = -\lambda z(t) + r(z(t)),
\end{equation*}
where $|r(z)|\le K|z|^2$ for a constant $K$ and all $z$.
Therefore,
\begin{align*}
\frac{\ln\delta}{\lambda}+t(\delta,q_+)&=\frac{\ln\delta}{\lambda}+t(\delta,S^{-t_0}q_+)+t_0\\
&=-\int_\delta^{z(t_0)}\frac{dy}{\lambda y}+\int_{\delta}^{z(t_0)}\frac{dy}{\lambda y+r(y)}+\frac{\ln z(t_0)}{\lambda}\\
&= -\int_\delta^{z(t_0)}\frac{r(y)dy}{\lambda y(\lambda y+r(y))}+\frac{\ln z(t_0)}{\lambda},
\end{align*}
where $t(\delta,S^{-t_0}q_+)$ is defined analogously to~\eqref{eq:definition_of_t_q_+}.
Our claim follows from the convergence of the integral in the r.h.s. as $\delta\to 0$.

Proof of Lemma~\ref{lm:growth_of_Phi}. Let us choose a new basis in $\R^d$ so that in that basis the
Euclidean norm of $A$ (denoted by $\|A\|$) is bounded by $\lambda+\alpha/2$. This
can be done as in \cite[Section 1.2]{Katok--Hasselblatt}. Now Lemma 4.1 from \cite[Chapter IV]{Hartman}
implies
\begin{equation*}
\|\Phi_s(t)\|\le e^{\int_s^t\|A(r)\|dr},
\end{equation*}
and our claim follows from $\lim_{r\to\infty}A(r)=A$.

Proof of Lemma~\ref{lm:growth_of_Y}. It\^o's isometry implies
\begin{equation*}
\E |Y(t)|^2=\E \left(\int_0^t\Phi_r(t)dW(r)\right)^2=\int_0^t \|\Phi_r(t)\|^2_2 dr,
\end{equation*}
where $\|B\|_2^2=\sum_{i,j} B_{ij}^2$ is the square of the quadratic norm of a matrix $B$, so that
due to Lemma~\ref{lm:growth_of_Phi} and the equivalence of any two norms in $\R^d$, we have 
\begin{equation}
\label{eq:growth_of_second_order_of_Y}
\E |Y(t)|^2\le K'_{\alpha/2} \int_0^t e^{2(\lambda+\alpha/2)(t-r)} dr\le K''_{\alpha/2} e^{2(\lambda+\alpha/2)t},
\end{equation}
for some constants $K'_{\alpha/2}$, $K''_{\alpha/2}$, and all $t\ge 0$. Inequality~\eqref{eq:growth_of_second_order_of_Y}
with the Borel--Cantelli Lemma implies the desired growth of $Y(t)$ along integer values of $t$. To interpolate between the integer times,
we apply the standard Kolmogorov--Chentsov technique based on an estimate for increments of $Y$. For any $z>0$, we have
\begin{align*}
\Pp\{|Y(t_2)-Y(t_1)|\ge z\}&\le \Pp\left\{\int_{t_1}^{t_2} A(s)Y(s)ds>\frac{z}{2}\right\}+\Pp\left\{|W(t_2)-W(t_1)|>\frac{z}{2}\right\}\\
&\le\frac{4}{z^2}\E \left(\int_{t_1}^{t_2} A(s)Y(s)ds\right)^2+\frac{16}{z^4}\E\left(W(t_2)-W(t_1)\right)^4\\
&\le\hat K_{\alpha/2} (t_2-t_1)^2 \left[\frac{1}{z^2} e^{2(\lambda+\alpha/2)t_2}+\frac{1}{z^4}\right],
\end{align*}
where $\hat K_{\alpha/2}$ is a positive constant.

For $n,m\in \{0\}\cup\N$, we introduce $D_{m,n}$ as the set of all the rationals of the form $k/2^n\in[m,m+1]$ with integer $k$. For each $t\in D_{m,n}$ with $n\in\N$,
we define $t_-=\sup\{s\in D_{m,n-1}: s\le t\}$. Then $|t-t_-|\le 2^{-n}$.
Pick any $\rho$ with $1<\rho^4<2$. The continuity of the trajectories of $Y$ implies that
\begin{multline*}
\Pp\left\{\sup_{s\in[m,m+1]}|Y(s)|\ge e^{(\lambda+\alpha)m}\right\} 
\\ \le \Pp\left\{|Y(m)|\ge e^{(\lambda+\alpha)m}\frac{\rho-1}{\rho}\right\}
                              +\sum_{n=1}^{\infty}\,\sum_{t\in D_{m,n}}\Pp\left\{|Y(t)-Y(t_-)|\ge e^{(\lambda+\alpha)m}\frac{\rho-1}{\rho^{n+1}}\right\}
\\ \le \frac{K''_{\alpha/2} e^{2(\lambda+\alpha/2)m}\rho^2}{e^{2(\lambda+\alpha)m}(\rho-1)^2}
+\sum_{n=1}^{\infty} 2^{n}\hat K_{\alpha/2} (2^{-n})^2 \left[\frac{e^{2(\lambda+\alpha/2)(m+1)}\rho^{2n+2}}{e^{2(\lambda+\alpha)m}(\rho-1)^2 } +\frac{\rho^{4n+4}}{e^{4(\lambda+\alpha)m}(\rho-1)^4}\right].
\end{multline*}
Due to our choice of $\rho$, the series in the r.h.s. converges exponentially, and the whole r.h.s. decays in $m$ as $e^{-\alpha m}$, so that
we can finish the proof applying the Borel--Cantelli lemma.

\bibliographystyle{alpha}
%\bibliography{happydle}
\def\cprime{$'$}

\end{document}